\documentclass[12pt]{article}
\usepackage{amsmath,amssymb,amsthm}
\allowdisplaybreaks
\def\cc{\mathbb{C}}
\def\M{\mathcal M}
\def\N{\mathcal N}
\newtheorem{set}{set}[section]
\newtheorem{Corollary}[set]{Corollary}

\newtheorem{Lemma}[set]{Lemma}
\newtheorem{Proposition}[set]{Proposition}
\newtheorem{Remark}[set]{Remark}
\newtheorem{Theorem}[set]{Theorem}
\newtheorem{Example}[set]{Example}
\numberwithin{equation}{section}
\begin{document}
\pagestyle{myheadings}
\date{}
\title{On spectra and Brown's spectral measures of elements in free products of matrix algebras}
\author{Junsheng Fang \and Don Hadwin
 \and Xiujuan Ma }
\maketitle
\begin{abstract}
We compute spectra and Brown measures of some non self-adjoint
operators in $(M_2(\cc), \frac{1}{2}Tr)*(M_2(\cc), \frac{1}{2}Tr)$,
the reduced free product von Neumann algebra of $M_2(\cc)$ with
$M_2(\cc)$. Examples include $AB$ and $A+B$, where $A$ and $B$ are
matrices in $(M_2(\cc), \frac{1}{2}Tr)*1$ and $1*(M_2(\cc),
\frac{1}{2}Tr)$, respectively. We prove that $AB$ is an R-diagonal
operator (in the sense of Nica and Speicher~\cite{N-S1}) if and only
if $Tr(A)=Tr(B)=0$.  We show that if $X=AB$ or  $X=A+B$ and $A,B$
are not scalar matrices, then the Brown measure of $X$ is not
concentrated on a single point. By a theorem of Haagerup and
Schultz~\cite{H-S1}, we obtain that if $X=AB$ or $X=A+B$ and $X\neq
\lambda 1$, then $X$ has a nontrivial hyperinvariant subspace
affiliated with $(M_2(\cc), \frac{1}{2}Tr)*(M_2(\cc),
\frac{1}{2}Tr)$.
\end{abstract}

\section{Introduction} In 1983, L.G.~Brown~\cite{Br} introduced a spectral
distribution measure for non-normal elements in a finite von Neumann
algebra with respect to a fixed normal faithful tracial state, which
is called the Brown measure of the operator.  Recently, U.~Haagerup
and H.~Schultz~\cite{H-S1} proved a remarkable result which states
that if the support of Brown measure of an operator in a type ${\rm
II}_1$ factor contains more than two points, then the operator has a
non-trivial hyperinvariant subspace affiliated with the type ${\rm
II}_1$ factor. In general cases, the computation of Brown measures
of non-normal operators  are nontrivial. The first essential result
was given by Haagerup and F.~Larsen. In~\cite{H-L}, Haagerup and
Larsen computed the spectrum and Brown measure of R-diagonal
operators in a finite von Neumann algebra, in terms of the
distribution of its radial part. Brown measures of some non-normal
and non-R-diagonal operators, examples include $u_n+u_\infty$, where
$u_n$ and $u_\infty$ are the generators of $\mathbb{Z}_n$ and
$\mathbb{Z}$ respectively, in the free product
$\mathbb{Z}_n*\mathbb{Z}$, and elements of the form
$S_\alpha+iS_\beta$, where $S_\alpha$ and $S_\beta$ are free
semi-circular elements of variance $\alpha$ and $\beta$, are
computed by P.~Biane and F.~Lehner in~\cite{B-L}. The purpose of
this paper is to compute the spectra and  Brown measures of some non
hermitian operators in $(M_2(\cc), \frac{1}{2}Tr)*(M_2(\cc),
\frac{1}{2}Tr)$, the reduced free product von Neumann algebra of
$M_2(\cc)$ with $M_2(\cc)$ (cf. [Ch]). Examples include $AB$ and
$A+B$, where $A$ and $B$ are matrices in $(M_2(\cc),
\frac{1}{2}Tr)*1$ and $1*(M_2(\cc),
\frac{1}{2}Tr)$, respectively. This paper is organized as follows.\\

In section 2 we recall preliminary facts about Brown measures,
R-diagonal operators, Haagerup and Larsen's result on Brown measures
of R-diagonal operators and some notation used in this paper. In
section 3, we provide some results on the spectra and spectral
radius of operators in $M_2(\cc)*M_2(\cc)$, the universal free
product C*-algebra of $M_2(\cc)$ with $M_2(\cc)$. Firstly we compute
the spectral radius of $AB$ for two normal matrices $A\in
M_2(\cc)*1$ and $B\in 1*M_2(\cc)$ relative to $M_2(\cc)*M_2(\cc)$.
As a corollary, we also get the spectrum radius of $AB$ for normal
matrices $A\in (M_2(\cc), \frac{1}{2}Tr)*1$ and $B\in 1*(M_2(\cc),
\frac{1}{2}Tr)$, relative to the reduced free product von Neumann
algebra of $M_2(\cc)$ with $M_2(\cc)$. Then we obtain the following
 result: Let $A$,$B$ be matrices in $M_2(\cc)*1$ and $1*M_2(\cc)$,
respectively, such that $Tr(A)=Tr(B)=0$. Then $\sigma(AB)$, the
spectrum of $AB$, relative to $M_2(\cc)*M_2(\cc)$, is the closure of
the annulus centered at 0 with inner radius
$\|A^{-1}\|^{-1}\|B^{-1}\|^{-1}$ and outer radius $\|A\|\|B\|$,
where we use the convention $\infty^{-1}=0$ and if $A$ is not
invertible then $\|A^{-1}\|:=\infty$.\\

In section 4 we prove that $AB$ is an R-diagonal operator if and
only if $Tr(A)=Tr(B)=0$, where $A\in (M_2(\cc), \frac{1}{2}Tr)*1$
and $B\in 1* (M_2(\cc), \frac{1}{2}Tr)$. As a corollary, we
explicitly compute the spectrum and  Brown measure of $AB$
($Tr(A)=Tr(B)=0$) in terms of
$\S-$transform of $A^*A$ and $B^*B$.\\

 In
section 5, we develop algebraic techniques used in~\cite{Dy}.  Let
$X\in 1*(M_2(\cc), \frac{1}{2}Tr)$. With respect to the matrix units
of $(M_2(\cc), \frac{1}{2}Tr)*1$,
$X=\left(\begin{array}{cc} x_1&x_2\\
x_3&x_4
\end{array}\right).$ By~\cite{Dy}, $(M_2(\cc),
\frac{1}{2}Tr)*(M_2(\cc), \frac{1}{2}Tr)\cong L(\mathbb{F}_3)\otimes
M_2(\cc)$. So $x_1,x_2,x_3,x_4\in L(\mathbb{F}_3)$. In section 5, we
find  $\ast$-free generators $h,u,v$ of $L(\mathbb{F}_3)$ (different
from the free generators given in~\cite{Dy}) so that we may
explicitly write out $x_1,x_2,x_3, x_4$ in terms of $h,u,v$.\\

In section 6, we compute miscellaneous examples of Brown measures of
operators $A+B$ and $AB$, where $A\in (M_2(\cc), \frac{1}{2}Tr)*1$
and $B\in 1* (M_2(\cc), \frac{1}{2}Tr)$. As a corollary, we show
that  $A+B$ is an R-diagonal operator if and
only if $A+B=0$.\\

In section 7, we prove the following result: Let $A\in (M_2(\cc),
\frac{1}{2}Tr)*1$ and $B\in 1* (M_2(\cc),\frac{1}{2}Tr)$.  if
$X=A+B$ or $X=AB$ and $A,B$ are not scalar matrices, then the Brown
measure of $X$ is not concentrated on a single point. As a corollary
of Theorem 7.1 of [H-S1], we prove that if $X=A+B$ or $X=AB$ and
$X\neq \lambda 1$, then $X$ has a nontrivial hyperinvariant subspace
affiliated with $(M_2(\cc), \frac{1}{2}Tr)*(M_2(\cc),
\frac{1}{2}Tr)$.\\

Many concrete examples of spectra and Brown measures are given in
this paper.
 For some interesting applications, we refer to~\cite{FDM}.\\

\section{Preliminaries}
\subsection{Fuglede-Kadison determinant and Brown's spectral
measure.} Let $\M$ be a finite von Neumann algebra with a faithful
tracial state $\tau$. The \emph{Fuglede-Kadison
determinant}~\cite{F-K}, $\Delta:\,\M\rightarrow [0,\infty[$, is
given by
\[\Delta(T)=\exp\{\tau(\log|T|)\},\qquad T\in\M,\]
with $\exp\{-\infty\}:=0$. For an arbitrary element $T$ in $\M$ the
function $\lambda\rightarrow \log \Delta(a-\lambda 1)$ is
subharmonic on $\mathbb{C}$, and it's Laplacian
\[d\mu_{T}(\lambda):=\frac{1}{2\pi} \bigtriangledown^2 \log
\Delta(T-\lambda 1),\] in the distribution sense, defines a
probability measure $\mu_T$ on $\mathbb{C}$, called the
\emph{Brown's measure}~\cite{Br} of $T$.  From the definition, Brown
measure $\mu_T$ only depends on the joint
 distribution of $T$ and $T^*$.\\

If $T$ is normal, $\mu_T$
 is the trace $\tau$ composed with the spectral projections of $T$.
 If $\M=M_n(\cc)$ and $\tau=\frac{1}{n} Tr$ is the normalized trace
on $M_n(\cc)$, then $\mu_T$ is the normalized counting measure
$\frac{1}{n}\left(\delta_{\lambda_1}+\delta_{\lambda_2}\cdots+\delta_{\lambda_n}\right)$,
 where $\lambda_1,\lambda_2 \cdots,\lambda_n$ are the eigenvalues of $T$
 repeated according to root multiplicity.
\\

The Brown measure has the following properties (see~\cite{Br,
H-S2}): $\mu_T$ is the unique compactly supported measure on $\cc$
such that $\log\Delta(T-\lambda 1)=\int_{\cc}
\log|z-\lambda|d\mu_T(z)$ for all $\lambda\in\cc$. The support of
$\mu_T$ is contained in $\sigma(T)$, the spectrum of $T$.
$\mu_{ST}=\mu_{TS}$ for arbitrary $S,T$ in $\M$, and if $f(z)$ is
analytic in a neighborhood of $\sigma(A)$, $\mu_{f(T)}=(\mu_T)_f$,
the push-forward measure of $\mu_T$ under the map
$\lambda\rightarrow f(\lambda)$. If $E\in\M$ is a projection such
that $E\in Lat T$, then with respect to $E, I-E$ we can write
\[T=\left(\begin{array}{cc}
A&B\\
0&C
\end{array}\right),\]
where $A=ETE$ and $C=(I-E)T(I-E)$ are elements of $\M_1=E\M E$ and
$\M_2=(I-E)\M (I-E)$, respectively. Let $\mu_{A}$ and $\mu_{C}$ be
the Brown measures of $A$ and $C$ computed relative to $\M_1$ and
$\M_2$, respectively. Then $\mu_T=\alpha\mu_A+(1-\alpha)\mu_C$,
where $\alpha=\tau(E)$. \\

For a generalization of Brown measures of sets of commuting
operators in a type ${\rm II}_1$ factor, we refer to~\cite{Sc}.

\subsection{R-diagonal operators}
In 1995, A.~Nica and S.~Speicher~\cite{N-S1}  introduced the class
of R-diagonal operators in non-commutative probability spaces.
Recall that an operator $T$ in a non-commutative probability space
is an $R-$diagonal operator if the $R-$transform $R_{\mu(T,T^*)}$ of
the joint distribution $\mu(T,T^*)$ of $T, T^*$ is of the form
\[R_{\mu(T,T^*)}(z_1,z_2)=\sum_{n=1}^{\infty}\alpha_n (z_1z_2)^n+
\sum_{n=1}^{\infty}\alpha_n (z_2z_1)^n.\] Nica and
Speicher~\cite{N-S1} proved that $T$ is an $R-$diagonal operator if
and only if $T$ has same $\ast$-distribution as product $UH$, where
$U$ and $H$ are $\ast$-free random variables in some tracial non
commutative probability space, $U$ is a Haar unitary operator and
$H$ is positive. If $T$ is an R-diagonal operator, then the
$\ast$-distribution of $T$ is uniquely determined by the
distribution of $T^*T=|T|^2$. If $T$ is an R-diagonal operator and
$S$ is $\ast$-free with $T$, then both $ST$ and $TS$ are R-diagonal
operators (see~\cite{N-S1}). If $T$ is an R-diagonal operator and
$n\in \mathbb{N}$, then $T^n$ is also an R-diagonal operator
(see~\cite{H-L,La}).  For other important properties of R-diagonal
operators, we refer to~\cite{H-L, La,N-S1, N-S2}.

\subsection{Brown measures of R-diagonal operators}

In~\cite{H-L}, Haagerup and Larson explicitly computed the Brown
measures of {R-diagonal} operators in a finite von Neumann algebra.
\begin{Theorem} \emph{(Theorem 4.4 of~\cite{H-L})} Let $U,H$ be $\ast$-free random variables in a noncommutative probability space $(\M,\tau)$,
with $U$ a Haar unitary operator and $H$ a positive operator such
that the distribution $\mu_H$ of $H$ is not a Dirac measure. Then
the Brown measure $\mu_{UH}$ of $UH$ can be computed as the
following.
\begin{enumerate}
\item $\mu_{UH}$ is rotation invariant and its support is the
annulus with inner radius $\|H^{-1}\|_2^{-1}$ and outer radius
$\|H\|_2$.
\item $\mu_{UH}(\{0\})=\mu_{H}(\{0\})$ and for $t\in ]\mu_H(\{0\}),
1]$,
\[\mu_{UH}\left(\mathbb{B}\left(0, \left(\S_{\mu_{H^2}}(t-1)\right)^{-1/2}\right)\right)=t,\]
where $\S_{\mu_{H^2}}$ is the $\S-$transform of $H^2$ and
$\mathbb{B}(0,r)$ is the open disc with center 0 and radius $r$;
\item $\mu_{UH}$ is the only rotation invariant symmetric
probability measure satisfying 2.
\end{enumerate}
Furthermore, if $H$ is invertible, then $\sigma(UH)=supp\mu_{UH}$;
if $H$ is not invertible, then
$\sigma(UH)=\overline{\mathbb{B}(0,\|H\|_2)}.$
\end{Theorem}

\subsection{Some Notation}

The following notation will be used in the rest of the paper
\begin{itemize}
\item $(\M,\tau)=(M_2(\cc), \frac{1}{2}Tr)* (M_2(\cc), \frac{1}{2}Tr)$ denotes the reduced free product von Neumann algebra
of $M_2(\cc)$ with $M_2(\cc)$ with the unique tracial state $\tau$;
\item $M_2(\cc)_{(1)}:= (M_2(\cc), \frac{1}{2}Tr)*1$ and $M_2(\cc)_{(2)}:=
1*(M_2(\cc), \frac{1}{2}Tr)$;
\item $\{E_{ij}\}_{i,j=1,2}, \{F_{ij}\}_{i,j=1,2}$ are matrix units of
$M_2(\cc)_{(1)}$ and $M_2(\cc)_{(2)}$, respectively;
\item $P=E_{11}$ and $Q=F_{11}$;
\item $\M\cong \N\otimes M_2(\cc)_{(1)}\cong \N\otimes M_2(\cc)_{(2)}$.
For $X\in \M$,
$X=\left(\begin{array}{cc} x_1&x_2\\
x_3&x_4
\end{array}\right)_{(1)}=\left(\begin{array}{cc} x_1'&x_2'\\
x_3'&x_4'
\end{array}\right)_{(2)}$ means the decomposition is
with respect to above matrix units of $M_2(\cc)_{(1)}$ and
$M_2(\cc)_{(2)}$, respectively.
\item $W_0=\left(\begin{array}{cc}
1&0\\
0&1
\end{array}\right)_{(1)}, W_1=\left(\begin{array}{cc}
1&0\\
0&-1
\end{array}\right)_{(1)}, W_2=\left(\begin{array}{cc}
0&-1\\
1&0
\end{array}\right)_{(1)}, W_3=\left(\begin{array}{cc}
0&1\\
1&0
\end{array}\right)_{(1)}$;
\item $V_0=\left(\begin{array}{cc}
1&0\\
0&1
\end{array}\right)_{(2)}, V_1=\left(\begin{array}{cc}
1&0\\
0&-1
\end{array}\right)_{(2)}, V_2=\left(\begin{array}{cc}
0&-1\\
1&0
\end{array}\right)_{(2)}, V_3=\left(\begin{array}{cc}
0&1\\
1&0
\end{array}\right)_{(2)}$;
\item $A,A_1,\cdots,A_n$ denote elements in $M_2(\cc)_{(1)}$,
$B,B_1,\cdots,B_n$ denote elements in $M_2(\cc)_{(2)}$, $X,Y,Z$
denote general elements in $\M$;
\item An element $X$ in $\M$ is called \emph{centered} if
$\tau(X)=0$.
\end{itemize}

We end this section with the following lemma. The proof is an easy
exercise.
\begin{Lemma}  $V_1M_2(\cc)_{(1)} V_1$ is free with
$M_2(\cc)_{(1)}$.
\end{Lemma}

\section{Spectra of elements in the universal free product of $M_2(\cc)$ and $M_2(\cc)$}
Let $\AA=M_2(\cc)*M_2(\cc)$ denote the universal free product
C*-algebra of $M_2(\cc)$ with $M_2(\cc)$. Then there is a *
homomorphism $\pi$ from  $\AA$ onto the reduced free product
C*-algebra of $M_2(\cc)$ and $M_2(\cc)$, the C*-subalgebra generated
by $M_2(\cc)_{(1)}$ and $M_2(\cc)_{(2)}$ in $\M$. Since
$\sigma(\pi(a))\subseteq \sigma(a)$ for $a\in \AA$, it is useful to
obtain some information of spectrum of $AB$, where $A\in M_2(\cc)*1$
and $B\in 1* M_2(\cc)$.

\subsection{``Free product" of normal matrices}
\begin{Lemma} Let $A\in M_2(\cc)*1$ and $B\in 1*M_2(\cc)$ be normal
matrices. Then $r(AB)=\|A\|\cdot \|B\|$ relative to $\AA$.
\end{Lemma}
\begin{proof} $r(AB)\leq \|AB\|\leq \|A\|\cdot\|B\|$.  We need only to prove $r(AB)\geq \|A\|\cdot\|B\|$.
Since $A$ is a normal matrix, there is a unitary matrix $U_1\in
M_2(\cc)*1$ such that $U_1AU_1^*=\left(\begin{array}{cc}
\alpha_1&0\\
0&\beta_1
\end{array}\right)$
 and $\|\alpha_1\|=\|A\|$. Similarly, there is a unitary matrix $U_2\in
1*M_2(\cc)$ such that $U_2BU_2^*=\left(\begin{array}{cc}
\alpha_2&0\\
0&\beta_2
\end{array}\right)$
 and $\|\alpha_2\|=\|B\|$. Let $\pi_1(X)=U_1XU_1^*$ and
 $\pi_2(Y)=U_2YU_2^*$ be $\ast$-representations of $M_1(\cc)*1$ and
 $1*M_2(\cc)$ to $M_2(\cc)$, respectively. Then there is a
 $\ast$-representation $\pi=\pi_1*\pi_2$ from $\AA$ to $M_2(\cc)$ and
 $\pi(AB)=\left(\begin{array}{cc}
 \alpha_1\alpha_2&0\\
 0&\beta_1\beta_2
 \end{array}\right)$.  Therefore, $\alpha_1\alpha_2\in
 \sigma(\pi(AB))\subseteq \sigma(AB)$. So $r(AB)\geq
 |\alpha_1\alpha_2|=\|A\|\cdot \|B\|$.
\end{proof}
\begin{Corollary}Let $A\in M_2(\cc)_{(1)}$ and $B\in M_2(\cc)_{(2)}$ be normal
matrices. Then $r(AB)=\|A\|\cdot \|B\|$ relative to $\M$.
\end{Corollary}
\begin{proof}   We may assume
that $A$ and $B$ are diagonal matrices. Then we can treat $AB$ as an
operator in the full free product $C^*(\mathbb{Z}_2\ast
\mathbb{Z}_2)$. Same technique used in the previous lemma gives the
corollary.
\end{proof}

\subsection{``Free product" of non-normal matrices}

It is well-known that two matrices $X,Y$ in $M_2(\cc)$ are unitarily
equivalent if and only if $Tr(X)=Tr(Y), Tr(X^2)=Tr(Y^2)$ and
$Tr(X^*X)=Tr(Y^*Y)$.  The proof of the following lemma now is an
easy exercise.
\begin{Lemma} If $A\in M_2(\cc)$ and $Tr(A)=0$, then
$A$ is unitarily equivalent to a matrix of form $\displaystyle
\left(\begin{array}{cc} 0&\alpha\\
\beta&0
\end{array}\right),$ where $\alpha,\beta$ are complex numbers.
\end{Lemma}

\begin{Remark} We have the following useful observations:
\begin{itemize}
\item $\left(\begin{array}{cc}
0&1\\
1&0
\end{array}\right)\left(\begin{array}{cc}
0&\alpha\\
\beta&0
\end{array}\right)\left(\begin{array}{cc}
0&1\\
1&0
\end{array}\right)=\left(\begin{array}{cc}
0&\beta\\
\alpha&0
\end{array}\right).$
\item $\left(\begin{array}{cc}
1&0\\
0&e^{i(\theta_1-\theta_2)/2}
\end{array}\right)\left(\begin{array}{cc}
0&|\alpha|e^{i\theta_1}\\
|\beta|e^{i\theta_2}&0
\end{array}\right)\left(\begin{array}{cc}
1&0\\
0&e^{-i(\theta_1-\theta_2)/2}
\end{array}\right)=e^{i(\theta_1+\theta_2)/2}\left(\begin{array}{cc}
0&|\alpha|\\
|\beta|&0
\end{array}\right).$
\end{itemize}
\end{Remark}

\begin{Lemma}Let $A\in M_2(\cc)*1$ and $B\in 1*M_2(\cc)$ be
matrices such that $Tr(A)=Tr(B)=0$. Then $r(AB)=\|A\|\cdot \|B\|$
relative to $\AA$.
\end{Lemma}
\begin{proof} We need only to prove $r(AB)\geq \|A\|\cdot \|B\|$.
 By Lemma 3.3 and Remark 3.4, there are unitary matrices
$U,V$ in $M_2(\cc)$ such that
$UAU^*=\left(\begin{array}{cc} 0&\alpha_1\\
\beta_1&0\\
\end{array}\right)$ and $VBV^*=\left(\begin{array}{cc} 0&\alpha_2\\
\beta_2&0\\
\end{array}\right)$ and $|\alpha_1|=\|A\|, |\beta_2|=\|B\|$. Let
$\pi_1(X)=UXU^*$ and $\pi_2(Y)=VYV^*$ be $\ast$-representations of
$M_2(\cc)*1$ and $1*M_2(\cc)$ to $M_2(\cc)$, respectively. Let
$\pi=\pi_1*\pi_2$ be the induced $\ast$-representation of $\AA$ to
$M_2(\cc)$. Then $\sigma(AB)\supseteq
\sigma(\pi(AB))=\sigma(\pi_1(A)\pi_2(B))=\{\alpha_1\beta_2,
\alpha_2\beta_1\}$. Therefore, $r(AB)\geq
|\alpha_1\beta_2|=\|A\|\cdot \|B\|$.
\end{proof}

\begin{Theorem}Let $A\in M_2(\cc)*1$ and $B\in 1*M_2(\cc)$ be
matrices such that $Tr(A)=Tr(B)=0$. Then
\[\sigma(AB)=[\|A^{-1}\|^{-1}\|B^{-1}\|^{-1},
\|A\|\|B\|]\times_p [0,2\pi],\] where  $\times_p$ denotes the polar
set product $\{re^{i\theta}:\,\, r\in
[\|A^{-1}\|^{-1}\|B^{-1}\|^{-1}, \|A\|\|B\|],\,\, \theta\in
[0,2\pi]\}$.
\end{Theorem}
\begin{proof} We will prove the theorem for two cases.\\

\noindent Case 1. Either $A$ or $B$ is not invertible. We may assume
that $A$ is not invertible. By  $Tr(A)=0$, Lemma 3.3 and Remark 3.4,
$A$ is unitarily
equivalent to $\left(\begin{array}{cc} 0&\alpha_1\\
0&0
\end{array}\right)$. Without loss of generality, we assume that $A=\left(\begin{array}{cc} 0&1\\
0&0
\end{array}\right)\in M_2(\cc)*1$. By Lemma 3.3 and Remark 3.4, we may also assume that $B=\left(\begin{array}{cc} 0&\alpha\\
\beta&0
\end{array}\right)\in 1*M_2(\cc)$ and $\beta\geq \alpha\geq 0$.  We need to prove
that $\sigma(AB)$ is the closed disc of complex plane with center 0
and radius $\beta$. Since $A$ is unitarily equivalent to
$e^{i\theta}A$ in $M_2(\cc)*1$, $\sigma(AB)$ is rotation invariant.
For $\theta\in [0,2\pi]$, let $U=\left(\begin{array}{cc} \cos\theta& \sin\theta\\
-\sin\theta&\cos\theta
\end{array}\right)$. Let $\pi_1(X)=X$ and $\pi_2(Y)=UYU^*$ be
$\ast$-representations of $M_2(\cc)*1$ and $1*M_2(\cc)$ to
$M_2(\cc)$, respectively. Let $\pi=\pi_1*\pi_2$ be the induced
$\ast$-representation of $\AA$ to $M_2(\cc)$. Then
\[\pi(AB)=AUBU^*=\left(\begin{array}{cc}
-\alpha
\sin^2\theta+\beta\cos^2\theta&-(\alpha+\beta)\sin\theta\cos\theta\\
0&0
\end{array}\right).\] So $\sigma(\pi(AB))=\{-\alpha
\sin^2\theta+\beta\cos^2\theta, 0\}$. Since $[0,\beta]\subseteq
[-\alpha,\beta]=\{-\alpha \sin^2\theta+\beta\cos^2\theta:\,
\theta\in [0,2\pi]\}$, $[0,\beta]\subseteq \sigma(AB)$. Since
$\sigma(AB)$ is rotation invariant, $\sigma(AB)$ contains the closed
disc with center 0 and radius $\beta$. By Lemma 3.5, $\sigma(AB)$ is
the closed disc of complex plane with center 0 and
radius $\beta$.\\

 \noindent Case 2. Both $A$ and $B$ are invertible. By Lemma 3.3 and
Lemma 3.4, we may assume that $A=\left(\begin{array}{cc} 0&1\\
\beta_1&0 \end{array}\right) $ and $B=\left(\begin{array}{cc}0&1\\
\beta_2&0\end{array}\right)$ such that $\beta_1,\beta_2\geq 1$. Then
$A^{-1}=\left(\begin{array}{cc} 0&\beta_1^{-1}\\
1&0 \end{array}\right)$ and $B^{-1}=\left(\begin{array}{cc} 0&\beta_2^{-1}\\
1&0 \end{array}\right)$. We need to prove that
$\sigma(AB)=[1,\beta_1\beta_2]\times_p [0,2\pi]$. By Lemma 3.5,
$r(AB)=\beta_1\beta_2$ and $r((AB)^{-1})=1$. This implies that
$\sigma(AB)\subseteq [1,\beta_1\beta_2]\times_p [0,2\pi]$.  So we
need only to prove $\sigma(AB)\supseteq [1,\beta_1\beta_2]\times_p
[0,2\pi]$.\\

For $\phi,\psi\in [0,2\pi]$, let $U=\left(\begin{array}{cc}
\cos\psi& e^{i\phi}\sin\psi\\
-\sin\psi& e^{i\phi}\cos\psi
\end{array}\right)$. Then $U$ is a unitary matrix. Let $\pi_1(X)=UXU^*$ and $\pi_2(Y)=Y$ be
$\ast$-representations of $M_2(\cc)*1$ and $1*M_2(\cc)$ to
$M_2(\cc)$, respectively. Let $\pi=\pi_1*\pi_2$ be the induced
$\ast$-representation of $\AA$ to $M_2(\cc)$. Then
\[\pi(AB)=\left(\begin{array}{cc}
-\beta_1\beta_2e^{i\phi}\sin^2\psi+\beta_2e^{-i\phi}\cos^2\psi&*\\
*&\beta_1e^{i\phi}\cos^2\psi-e^{-i\phi}\sin^2\psi
\end{array}\right).\]
Let $\lambda_1(\phi,\psi),\lambda_2(\phi,\psi)$ be the eigenvalues
of $\pi(AB)$. Then
\begin{equation}
\lambda_1(\phi,\psi)\lambda_2(\phi,\psi)=\det
(\pi(AB))=\det(A)\det(B)=\beta_1\beta_2,
\end{equation}
\begin{equation}
\lambda_1(\phi,\psi)+\lambda_2(\phi,\psi)=(\beta_1e^{i\phi}+\beta_2e^{-i\phi})\cos^2\psi-
(\beta_1\beta_2e^{i\phi}+e^{-i\phi})\sin^2\psi.
\end{equation}
Note that $\sigma(AB)\supseteq \{\lambda_1(\phi,\psi):\,
\phi,\psi\in [0,2\pi]\}$.  We only need to prove that
$\{\lambda_1(\phi,\psi):\, \phi,\psi\in [0,2\pi]\}\supseteq
[1,\beta_1\beta_2]\times_p [0,2\pi]$. For this purpose, we need to
show for any $r\in [1,\beta_1\beta_2]$, $\theta\in [0,2\pi]$, there
are $\phi,\psi\in [0,2\pi]$ such that
\begin{equation}
re^{i\theta}+\frac{\beta_1\beta_2}{r}e^{-i\theta}=
(\beta_1e^{i\phi}+\beta_2e^{-i\phi})\cos^2\psi-
(\beta_1\beta_2e^{i\phi}+e^{-i\phi})\sin^2\psi.
\end{equation}
Let $\alpha=\cos^2\psi$. Simple computations show that equation 3.3
is equivalent to the following
\[
\left(r+\frac{\beta_1\beta_2}{r}\right)\cos\theta+i\left(r-\frac{\beta_1\beta_2}{r}\right)\sin\theta=\]\[
(\alpha(1+\beta_1)(1+\beta_2)-(1+\beta_1\beta_2))\cos\phi+i(\alpha(\beta_1-1)(\beta_2+1)+
(1-\beta_1\beta_2))\sin\phi.\]

Let
\[\Omega_1=\left\{\left(r+\frac{\beta_1\beta_2}{r}\right)\cos\theta+
i\left(r-\frac{\beta_1\beta_2}{r}\right)\sin\theta: \,\, r\in
[1,\beta_1\beta_2],\,\theta\in [0,2\pi]\right\},\]
\[\Omega_2=\{(\alpha(1+\beta_1)(1+\beta_2)-(1+\beta_1\beta_2))\cos\phi+i(\alpha(\beta_1-1)(\beta_2+1)+
(1-\beta_1\beta_2))\sin\phi\]\[:\, \alpha\in [0,1],\,
\phi\in[0,2\pi]\}.\] Now we need only to prove $\Omega_1=\Omega_2$.
Note that $\Omega_1$ is the union of a family of ellipses with
center the origin point and semimajor axis and semiminor axis
$|r+\frac{\beta_1\beta_2}{r}|$ and $|r-\frac{\beta_1\beta_2}{r}|$,
$1\leq r\leq \beta_1\beta_2$, respectively. Similarly, $\Omega_2$ is
the union of a family of ellipses with center the origin point and
semimajor axis and semiminor axis
$|\alpha(1+\beta_1)(1+\beta_2)-(1+\beta_1\beta_2)|$ and
$|\alpha(\beta_1-1)(\beta_2+1)+ (1-\beta_1\beta_2)|$, $0\leq
\alpha\leq 1$,  respectively. Note that the ``largest" ellipse in
$\Omega_1$ is with semimajor axis and semiminor axis
$|1+\beta_1\beta_2|$ and $|\beta_1\beta_2-1|$, respectively; the
``smallest" ellipse in $\Omega_1$ is with semimajor axis and
semiminor axis $2\sqrt{|\beta_1\beta_2|}$ and $0$, respectively. The
``largest" ellipse in $\Omega_2$ is with semimajor axis and
semiminor axis $|1+\beta_1\beta_2|$ and $|\beta_1\beta_2-1|$,
respectively; the ``smallest" ellipse in $\Omega_2$ is with
semimajor axis and semiminor axis $0$ and
$\frac{2\beta_1(\beta_2-1)}{\beta_1+1}$. So both $\Omega_1$ and
$\Omega_2$ are the closure of the domain enclosed by the ellipse
with center the origin point and semimajor axis and semiminor axis
$|1+\beta_1\beta_2|$ and $|\beta_1\beta_2-1|$, respectively. Thus
$\Omega_1=\Omega_2$.
\end{proof}

\section{R-diagonal operators in
$\M$}

In this section, we prove the following result. We will use the
notation introduced in section 2.4.

\begin{Theorem} In $\M$, let $A\in M_2(\cc)_{(1)}$ and $B\in M_2(\cc)_{(2)}$. Then
$AB$ is an R-diagonal operator if and only if $\tau(A)=\tau(B)=0$.
\end{Theorem}

To prove Theorem 4.1, we need the following lemmas.

\begin{Lemma} $\{W_1,V_1,W_3V_3\}''\cong L(\mathbb{Z}_2)*L(\mathbb{Z}_2)*L(\mathbb{Z}).$
\end{Lemma}
\begin{proof}  Let $U=W_3V_3$. Then $U$ is a Haar unitary operator.
We need only to prove that $U$ is * free with the von Neumann
subalgebra generated by $W_1$ and $V_1$. Let $g_1g_2\cdots g_n$ be
an alternating product of $\{U^n: n\neq 0\}$ and $\{W_1,V_1,
W_1V_1,V_1W_1, W_1V_1W_1,V_1W_1V_1,\cdots\}$. By regrouping, it is
an alternating product of $\{W_1, W_1W_3, W_3^*W_1, W_3^*W_1W_3,
W_3, W_3^*\}$ and $\{V_1, V_3V_1, V_1V_3^*, V_3V_1V_3^*, V_3,
V_3^*\}$. Thus the trace is 0.
\end{proof}

\begin{Lemma} $\displaystyle
\left(\begin{array}{cc} 0&\alpha_1\\
\beta_1&0
\end{array}\right)_{(1)}
\left(\begin{array}{cc} 0&\alpha_2\\
\beta_2&0
\end{array}\right)_{(2)}$ is an R-diagonal operator.
\end{Lemma}
\begin{proof} Note that
\[\left(\begin{array}{cc} 0&\alpha_1\\
\beta_1&0
\end{array}\right)_{(1)}
\left(\begin{array}{cc} 0&\alpha_2\\
\beta_2&0
\end{array}\right)_{(2)}=\left(\begin{array}{cc} \alpha_1&0\\
0&\beta_1
\end{array}\right)_{(1)}
\left(\begin{array}{cc} 0&1\\
1&0
\end{array}\right)_{(1)}\left(\begin{array}{cc} 0&1\\
1&0
\end{array}\right)_{(2)}
\left(\begin{array}{cc} \beta_2&0\\
0&\alpha_2
\end{array}\right)_{(2)}.
\]  By Lemma 4.2 and basic properties of R-diagonal operators given in 2.2, we prove the lemma.
\end{proof}

\begin{Lemma} With the assumption of Theorem 4.1 and assume $AB$ is an R-diagonal operator and
$\tau(A^2)\neq 0$. Then $\tau(B)=0$.
\end{Lemma}
\begin{proof}  Since $AB$ is an $R-$diagonal operator, $\tau(AB)=0$.
Since $A,B$ are $\ast$-free, $\tau(A)\tau(B)=\tau(AB)=0$. If
$\tau(B)=0$, then done. Otherwise, assume $\tau(A)=0$. Then
$0=\tau(ABAB)=\tau(A^2B)\tau(B)=\tau(A^2)(\tau(B))^2$. By
assumption, $\tau(B)=0$.
\end{proof}

\begin{Lemma} Let $B\in M_2(\cc)_{(2)}$ and $\lambda$ be any complex
number. Then $\sigma(E_{12}B)=\sigma(E_{12}(\lambda+B))$.
\end{Lemma}
\begin{proof} By Jacobson's theorem, \[\sigma(E_{12}(\lambda+B))\cup
\{0\}=\sigma (E_{11}E_{12}(\lambda+B))\cup\{0\}=\sigma
(E_{12}(\lambda+B)E_{11})\cup\{0\}\]\[=\sigma(E_{12}BE_{11})\cup\{0\}=\sigma(BE_{12})\cup\{0\}.\]
\end{proof}

\noindent\emph{Proof of Theorem 4.1.}\,\, If $\tau(A)=\tau(B)=0$,
then by Lemma 3.3 and Lemma 4.3, $AB$ is an R-diagonal operator.
Conversely, assume that $AB$ is an R-diagonal operator. Then
$0=\tau(AB)=\tau(A)\cdot\tau(B)$. So either $\tau(A)=0$ or
$\tau(B)=0$. Without loss of generality, we assume that $\tau(A)=0$.
If $\tau(A^2)\neq 0$, then $\tau(B)=0$ by Lemma 4.4. If
$\tau(A^2)=0$, then $A$ is unitary equivalent to $\alpha E_{12}$. We
may assume that $A=E_{12}$. By Theorem 2.1, if $E_{12}B$ is an
R-diagonal operator, then
$(r(E_{12}B))^2=\tau(B^*E_{21}E_{12}B)=\tau(E_{21}E_{12}BB^*)=\|E_{12}\|_2^2\cdot
\|B\|_2^2$. Note that $E_{12}(B-\tau(B))$ is an R-diagonal operator,
$(r(E_{12}(B-\tau(B)))^2=\|E_{12}\|_2^2\cdot \|B-\tau(B)\|_2^2$. By
Lemma 4.5, $\|B\|_2^2=\|B-\tau(B)\|_2^2$. This implies that
$\tau(B)=0$. This ends the proof.\\

Combining Theorem 4.1, Theorem 2.1 and  $\S$ transform of Voiculescu
(see~\cite{Vo, VDN}),  we have  the following theorem (It is
interesting to compare the following theorem and  Theorem 3.6).

\begin{Theorem} Let $A\in M_2(\cc)_{(1)}, B\in M_2(\cc)_{(2)}$ and $\tau(A)=\tau(B)=0$. Then
\begin{enumerate}
\item $\mu_{AB}$ is rotation invariant;
\item  $\sigma(AB)=supp\mu_{AB}=[\|A^{-1}\|_2^{-1}\|B^{-1}\|_2^{-1},
\|A\|_2\|B\|_2]\times_p [0,2\pi]$;
\item $\mu_{AB}(\{0\})=\max\{\mu_{A^*A}(\{0\}),\,\,\mu_{B^*B}(\{0\})\}$ and
\[\mu_{AB}(\mathbb{B}(0,
(\S_{\mu_{A^*A}}\S_{\mu_{B^*B}}(t-1))^{-1/2}))=t,\qquad
\text{for}\,\, t\in [\mu_{AB}(\{0\}), 1].\]
\end{enumerate}
\end{Theorem}

\section{Algebraic techniques}
 For $X\in \M$, define
\[\Phi(X)=\left(\begin{array}{cc}
E_{11}XE_{11}& E_{11}XE_{21}\\
E_{12}XE_{11}& E_{12}XE_{21}
\end{array}\right).\]
Then $\Phi$ is a $\ast$-isomorphism from $\M$ onto $E_{11}\M
E_{11}\otimes M_2(\cc)_{(1)}.$ We will identify $\M$ with $E_{11}\M
E_{11}\otimes M_2(\cc)_{(1)}$ by the canonical isomorphism $\Phi$.
In~\cite{Dy}, K.~Dykema proved that $E_{11}\M E_{11}\cong
L(\mathbb{F}_3)$. For $B\in M_2(\cc)_{(2)}$, we may write
\[B=\left(\begin{array}{cc}
b_{11}&b_{12}\\
b_{21}&b_{22}
\end{array}\right)\]
with respect to matrix units in $\M_1$. Then $b_{ij}\in
L(\mathbb{F}_3)$. In this section, we will develop the algebraic
techniques used in~\cite{Dy}. Combining the matrix techniques, we
may explicitly express $b_{ij}$ in terms of free generators of
$L(\mathbb{F}_3)$. \\

  Let $\Lambda\{W_1,V_1\}$ be the set of
words generated by $W_1,V_1$. Note that $W_1^2=V_1^2=1$ and
$\tau(W_1)=\tau(V_1)=0$. The following observation is crucial
in~\cite{Dy}. The proof is an easy exercise.

\begin{Lemma}$\tau(g_1g_2\cdots g_n)=0$ for an alternating product
of  $\Lambda\{W_1,V_1\}\setminus \{1,W_1\}$ and $\{E_{12},E_{21}\}$.
\end{Lemma}

Recall that  $P=E_{11}$ and $Q=F_{11}$. Let $W$ be the ``polar" part
of $(1-P)QP$ and $U=E_{12}W$. The following corollary is a special
case of Theorem 3.5 of~\cite{Dy}.
\begin{Corollary}\label{C:dykema's corollary} $U$ is a Haar unitary operator in $\M_P=P\M P$
and $U$, $PQP$ are $\ast$-free in $\M_P$.
\end{Corollary}

With the canonical identification of $\M$ with $\M_P\otimes
M_2(\cc)_{(1)}$,
\[Q=\left(\begin{array}{cc}
PQP& \sqrt{PQP-(PQP)^2}U\\
U^* \sqrt{PQP-(PQP)^2}& U^*(1-PQP)U
\end{array}\right).\]
By~\cite{Vo}, the distribution of $PQP$ (relative to $\M_P$) is
non-atomic and the density function is
\begin{equation}
\displaystyle \rho (t)=\frac{1}{\pi }\frac{1}{\sqrt{\frac{1}{4}-(\frac{1}{2}%
-t)^{2}}}, \qquad 0\leq t\leq 1.
\end{equation}
 By Corollary~\ref{C:dykema's corollary}, the von Neumann
subalgebra $\M_1$ generated by $M_2(\cc)_{(1)}$ and $Q$ is
$\ast$-isomorphic to $L(\mathbb{F}_2)\otimes M_2(\cc)_{(1)}$. Since
$\M_1$ is also $\ast$-isomorphic to $M_2(\cc)* L(\mathbb{Z}_2)$,
$M_2(\cc)* L(\mathbb{Z}_2)\cong L(\mathbb{F}_2)\otimes M_2(\cc)$,
which is proved by Dykema in~\cite{Dy}.\\

 Since $V_1=2Q-1$,
\[V_1=\left(\begin{array}{cc}
2PQP-1& 2\sqrt{PQP-(PQP)^2}U\\
2U^* \sqrt{PQP-(PQP)^2}& U^*(1-2PQP)U
\end{array}\right).\]
Simple computation shows that the density function of $2PQP-1$ is
\[\displaystyle \sigma(t)=\frac{1}{\pi}\frac{1}{\sqrt{1-t^2}},\qquad -1\leq t\leq 1.\] Let $H=2PQP-1$, then
\[V_1=\left(\begin{array}{cc}
H&\sqrt{1-H^2} U\\
U^*\sqrt{1-H^2}&  -U^*HU
\end{array}\right).\]

 Let $H=V|H|$ be the polar decomposition of $H$. Since
$H$ is a symmetric selfadjoint operator, $V^2=1$ and $V$ is
independent with the von Neumann algebra generated by $|H|$ in the
classical probability sense. Let $h=|H|, u=VU, v=UV$. Then $u,v$ are
Haar unitary operators and the distribution of $h$ relative to
$\M_P$ is non-atomic.

\begin{Lemma} $h, u, v$ are * free.
\end{Lemma}
\begin{proof}Let $g_1g_2\cdots g_n$ be an alternating product of
elements of $\mathfrak{S}=\{|H|\}''\ominus \cc I$, $\{(VU)^n: n\neq
0\}$, $\{(UV)^n: n\neq 0\}$. By regrouping, it is an alternating
product of elements of $\{\mathfrak{S}, V, V\mathfrak{S},
\mathfrak{S}V, V\mathfrak{S}V\}$ and $\{U^n: n\neq 0\}$. Since $H$
and $U$ are $\ast$-free, $\{\mathfrak{S}, V, V\mathfrak{S},
\mathfrak{S}V, V\mathfrak{S}V\}$ and $\{U^n: n\neq 0\}$ are free.
Since $V$ and $\mathfrak{S}$ are independent, $\tau(VS)=\tau(SV)=0$
for $S\in \mathfrak{S}$. This implies that $\tau(g_1g_2\cdots
g_n)=0$.

\end{proof}

 By simple computations, we have the following.
\begin{equation}
V_1E_{11}V_1=\left(\begin{array}{cc}
h^2&\sqrt{1-h^2}h u\\
u^*h\sqrt{1-h^2}&  u^*(1-h^2)u
\end{array}\right), \label{E:ve11v}
\end{equation}

\begin{equation}
V_1E_{12}V_1=\left(\begin{array}{cc}
HU^*\sqrt{1-H^2}&-HU^*HU\\
U^*\sqrt{1-H^2}U^*\sqrt{1-H^2}&  -U^*\sqrt{1-H^2}U^*HU
\end{array}\right)\label{E:ve12v-}\end{equation}
\begin{equation}=
\left(\begin{array}{cc}
hv^*\sqrt{1-h^2}&-hv^*hu\\
u^*\sqrt{1-h^2}v^*\sqrt{1-h^2}&  -u^*\sqrt{1-h^2}v^*hu
\end{array}\right).\label{E:ve12v}
\end{equation}

By Lemma 2.2,  $\M\cong M_2(\cc)_{(1)}*(V_1M_2(\cc)_{(1)}V_1)\cong
\M_P\otimes M_2(\cc)_{(1)}$. With this isomorphism, $\M_P$ is the
von Neumann algebra generated by $h, u $ and $v$ by~(\ref{E:ve11v})
and~(\ref{E:ve12v}). So $\M_P\cong L(\mathbb{F}_3)$. By simple
computations, we have
\[V_1\left(\begin{array}{cc}
\alpha&\beta\\
\gamma&\sigma
\end{array}\right)_{(1)}V_1=\left(\begin{array}{cc}
b_{11}&b_{12}\\
b_{21}& b_{22}
\end{array}\right),\]
where
\[
\begin{array}{cl}
  b_{11}= & \sigma+(\alpha-\sigma)h^2+\gamma\sqrt{1-h^2}vh+\beta hv^*\sqrt{1-h^2}, \\
  b_{12}= & (\alpha-\sigma)h\sqrt{1-h^2}u+\gamma \sqrt{1-h^2}v\sqrt{1-h^2}u-\beta hv^*hu, \\
  b_{21}= & (\alpha-\sigma)u^*h\sqrt{1-h^2}-\gamma u^*hvh+\beta u^*\sqrt{1-h^2}
v^*\sqrt{1-h^2}, \\
  b_{22}=& \alpha+(\sigma-\alpha) u^*h^2u-\gamma
u^*hv\sqrt{1-h^2}u-\beta u^*\sqrt{1-h^2}v^*hu.
\end{array}
\]

\begin{Theorem}\label{T:algebraic techinque}$\M{\cong} L(\mathbb{F}_3)\otimes
M_2(\cc)_{(1)}$; furthermore, let $\displaystyle
B=\left(\begin{array}{cc}
\alpha&\beta\\
\gamma&\sigma
\end{array}\right)_{(2)}$ in $M_2(\cc)_{(2)}$, then with respect to the matrix units
$\{E_{ij}\}_{i,j=1,2}\subset M_2(\cc)_{(1)}$, $\displaystyle
B=\left(\begin{array}{cc}
b_{11}&b_{12}\\
b_{21}& b_{22}
\end{array}\right)_{(1)}$, where $b_{ij}$ are given as above.
\end{Theorem}

\begin{Example}\label{E:normal matrices}\emph{ In Theorem~\ref{T:algebraic techinque}, let $\beta=\gamma=0$. Then we
have
\[ \left(\begin{array}{cc}
\alpha&0\\
0&\sigma
\end{array}\right)_{(2)}=\left(\begin{array}{cc}
\sigma+(\alpha-\sigma)h^2&(\alpha-\sigma)h\sqrt{1-h^2}u\\
(\alpha-\sigma)u^*h\sqrt{1-h^2}&\alpha+(\sigma-\alpha)u^*h^2u
\end{array}\right)_{(1)}.\]}
\end{Example}

\begin{Example}\label{E:algebra techinque}\emph{ In Theorem~\ref{T:algebraic techinque}, let $\alpha=\sigma$ and $\gamma=0$. Then we
have
\[\left(\begin{array}{cc}
\alpha&\beta\\
0&\alpha
\end{array}\right)_{(2)}=\left(\begin{array}{cc}
\alpha+\beta hv^*\sqrt{1-h^2}&-\beta hv^*hu\\
\beta u^*\sqrt{1-h^2} v^*\sqrt{1-h^2}&\alpha-\beta
u^*\sqrt{1-h^2}v^*hu
\end{array}\right)_{(1)}.
\]}
\end{Example}

\begin{Remark}\emph{ By equation (\ref{E:ve11v}), the distribution of $h^2$ is the
distribution of $E_{11}V_1E_{11}V_1E_{11}$ relative to $M_P$. So the
distribution of $h^2$ is same as the distribution of $PQP$ (relative
to $\M_P$). By [Vo], the distribution of $PQP$ (relative to $\M_P$)
is non-atomic and the density
function is \[\displaystyle \rho (t)=\frac{1}{\pi }\frac{1}{\sqrt{\frac{1}{4}-(\frac{1}{2}%
-t)^{2}}}, \qquad 0\leq t\leq 1.\]}
\end{Remark}

\section{Miscellaneous examples}

\begin{Example}\label{E:E12+F12}\emph{ We compute the Brown spectrum of $\alpha
E_{12}+\beta F_{12}$. Let $F_{12}=\left(\begin{array}{cc}0&1\\
0&0
\end{array}\right)_{(2)}=\left(\begin{array}{cc}b_1&b_2\\
b_3&b_4
\end{array}\right)_{(1)}.$ Then
\[(\alpha E_{12}+\beta
F_{12})^2=\alpha\beta(E_{12}F_{12}+F_{12}E_{12})=\alpha\beta(E_{12}+F_{12})^2=
\alpha\beta\left(\begin{array}{cc} b_3&b_1+b_4\\
0&b_3\end{array}\right)_{(1)}.
\] So $\mu_{(\alpha E_{12}+\beta
F_{12})^2}=\mu_{\alpha\beta b_3}$. By equation~(\ref{E:ve12v-}), the
distribution of $b_3$ is same as the distribution of
$(U^*\sqrt{1-H^2})^2$. Since $U^*\sqrt{1-H^2}$ is an R-diagonal
operator, $(U^*\sqrt{1-H^2})^2$ is also an R-diagonal operator.
Since the distribution of $\alpha E_{12}+\beta F_{12}$ is rotation
invariant, $\mu_{\alpha E_{12}+\beta
F_{12}}=\mu_{\sqrt{|\alpha\beta|}b},$ where $b=U^*\sqrt{1-H^2}$.
Simple computations show that (or by Proposition 5.10 and Corollary
5.11 of~\cite{FDM})
\[d\mu_{b}(z)=\frac{1}{\pi }\frac{1}{(1-r^2)^2}drd\theta\qquad
0\leq r\leq \frac{1}{\sqrt{2}}.
\] Hence
\[d\mu_{\alpha E_{12}+\beta
F_{12}}(z)=d\mu_{\sqrt{|\alpha\beta|}b}(z)=\frac{1}{\pi
}\frac{|\alpha\beta|^{3/2}}{(|\alpha\beta|-r)^2}drd\theta\qquad
0\leq r\leq \sqrt{\frac{|\alpha\beta|}{2}}
\] and
\[\sigma(\alpha E_{12}+\beta
F_{12})=\overline{\mathbb{B}\left(0,\sqrt{\frac{|\alpha\beta|}{2}}\right)}.
\]}
\end{Example}

\begin{Corollary}\label{C:E12+F12}
$r(\alpha E_{12}+\beta F_{12})=\sqrt{\frac{|\alpha\beta|}{2}}.$
\end{Corollary}

\begin{Corollary} Let $A\in M_2(\cc)_{(1)}$ and $B\in M_2(\cc)_{(2)}$.  Then $A+B$ is
an R-diagonal operator if and only if $A+B=0$.
\end{Corollary}
\begin{proof} Indeed, if $A+B$ is an R-diagonal operator, then $\tau(A+B)=0$. So we may assume
that $\tau(A)=\tau(B)=0$. Let $\lambda, -\lambda$ and $\eta,-\eta$
be the spectra of $A$ and $B$, respectively. Then
$0=\tau((A+B)^2)=\tau(A^2)+\tau(B^2)=\lambda^2+\eta^2$. By simple
computation we have
$0=\tau((A+B)^4)=\tau(A^4)+\tau(B^4)+4\tau(A^2)\tau(B^2)=\lambda^4+\eta^4+4\lambda^2\eta^2=(\lambda^2+\eta^2)^2+2\lambda^2\eta^2=2\lambda^2\eta^2$.
Thus $\lambda=\eta=0$. This implies that $A$ and $B$ are unitary
equivalent to $\alpha E_{12}$ and $ \beta F_{12}$, respectively. By
 Corollary~\ref{C:E12+F12}, $r(A+B)=\sqrt{\frac{|\alpha\beta|}{2}}.$ On the other
hand, since we assume that $A+B$ is an R-diagonal operator, by
Theorem 2.1, $(r(A+B))^2=\tau((A^*+B^*)(A+B))=\tau((\bar{\alpha}
E_{21}+\bar{\beta} F_{21})(\alpha E_{12}+\beta
F_{12}))=\frac{|\alpha|^2+|\beta|^2}{2}$. So
$|\alpha|^2+|\beta|^2=|\alpha\beta|$. This implies that $\alpha=0$
and $\beta=0$. Hence $A+B=0$.
\end{proof}

\begin{Example}\label{E:P times ab}\emph{ We compute the spectrum and Brown spectrum of
\[ X=\left(\begin{array}{cc}
1&0\\
0&0
\end{array}\right)_{(1)}\left(\begin{array}{cc}
\alpha&\beta\\
0&\alpha
\end{array}\right)_{(2)}.\] By Example~\ref{E:algebra techinque}, we have the following
\[X=\displaystyle \left(\begin{array}{cc}
1&0\\
0&0
\end{array}\right)_{(1)}\left(\begin{array}{cc}
\alpha&\beta\\
0&\alpha
\end{array}\right)_{(2)}=\displaystyle \left(\begin{array}{cc}
1&0\\
0&0
\end{array}\right)_{(1)}\left(\begin{array}{cc}
\alpha+\beta hv^*\sqrt{1-h^2}&-\beta hv^*hu\\
\beta u^*\sqrt{1-h^2} v\sqrt{1-h^2}&\alpha-\beta
u^*\sqrt{1-h^2}v^*hu
\end{array}\right)_{(1)}\]\[=\left(\begin{array}{cc}
\alpha+\beta hv^*\sqrt{1-h^2}&-\beta hv^*hu\\
0&0
\end{array}\right)_{(1)}.\]
So $\sigma(X)=\{0\}\cup \sigma(\alpha+\beta hv^*\sqrt{1-h^2})$ and
$\mu_X=\frac{1}{2}\delta_0+\frac{1}{2}\mu_{\alpha+\beta
hv^*\sqrt{1-h^2}}$.  Note that
$\mu_{hv^*\sqrt{1-h^2}}=\mu_{v^*\sqrt{1-h^2}h}$ and
$v^*\sqrt{1-h^2}h$ is an R-diagonal operator. We have the following
computations:
\[\|\sqrt{1-h^2}h\|_2^2=\tau_P((1-h^2)h^2)=\tau_P((1-PQP)PQP)=\int_{0}^1
\frac{1}{\pi }\frac{t(1-t)dt}{\sqrt{\frac{1}{4}-(\frac{1}{2}%
-t)^{2}}}=\frac{1}{8},\]
\[\|(\sqrt{1-h^2}h)^{-1}\|_2^2=\tau_P(((1-h^2)h^2)^{-1})=\tau_P(((1-PQP)PQP)^{-1})\]\[=\int_{0}^1
\frac{1}{\pi }\frac{dt}{t(1-t)\sqrt{\frac{1}{4}-(\frac{1}{2}%
-t)^{2}}}=\infty.\]
 By Theorem 2.1,
\[\sigma(X)=supp\mu_X=\{0\}\cup \overline
{\mathbb{B}(\alpha,|\beta|/2\sqrt{2})}.\]}
\end{Example}

\begin{Example}\label{E:P times normal}\emph{ We compute the spectrum and Brown spectrum of
\[\displaystyle Y=\left(\begin{array}{cc}
0&1\\
0&0
\end{array}\right)_{(1)}\left(\begin{array}{cc}
\alpha&0\\
0&\beta
\end{array}\right)_{(2)}.\] By Example~\ref{E:normal matrices}, we have the following
\[Y=\left(\begin{array}{cc}
0&1\\
0&0
\end{array}\right)_{(1)}\left(\begin{array}{cc}
\alpha&0\\
0&\beta
\end{array}\right)_{(2)}=\left(\begin{array}{cc}
0&1\\
0&0
\end{array}\right)_{(1)}\left(\begin{array}{cc}
\beta+(\alpha-\beta)h^2&(\alpha-\beta)h\sqrt{1-h^2}u\\
(\alpha-\beta)u^*h\sqrt{1-h^2}&\alpha+(\beta-\alpha)u^*h^2u
\end{array}\right)_{(1)}\]\[=\left(\begin{array}{cc}
(\alpha-\beta)u^*h\sqrt{1-h^2}&\alpha+(\beta-\alpha)u^*h^2u\\
0&0
\end{array}\right)_{(1)}.\]
Since $u^*h\sqrt{1-h^2}$ is an R-diagonal operator, similar
computations as Example~\ref{E:P times ab}, we have
\[\sigma(Y)=supp \mu_Y=\overline
{\mathbb{B}(0,|\alpha-\beta|/2\sqrt{2})}.\]}
\end{Example}

\begin{Example}\label{E:1+E12 times 1+F12}\emph{We compute the spectrum and Brown spectrum of
\[\displaystyle Z=(1+\alpha E_{12})(1+\beta F_{12})=\left(\begin{array}{cc}
1&\alpha\\
0&1
\end{array}\right)_{(1)}\left(\begin{array}{cc}
1&\beta\\
0&1
\end{array}\right)_{(2)}.\] For $\lambda\in\cc$, we have
\[Z-\lambda 1=(1+\alpha E_{12})(1+\beta F_{12})-
\lambda (1+\alpha E_{12})(1-\alpha E_{12})=(1+\alpha
E_{12})(\lambda\alpha E_{12}+\beta F_{12}-(\lambda-1)).\] This
implies that $\lambda \in \sigma(Z)$ if and only if $\lambda -1\in
\sigma (\lambda\alpha E_{12}+\beta F_{12})$. By
Example~\ref{E:E12+F12}, $\lambda -1\in \sigma (\lambda\alpha
E_{12}+\beta F_{12})$ if and only if
\[|\lambda-1|^2\leq \frac{ |\alpha\beta||\lambda|}{2}.\] So
\[\sigma(Z)=\left\{\lambda\in \cc:\,\, |\lambda-1|^2\leq \frac{
|\alpha\beta||\lambda|}{2}\right\}.\] In the following, we will show
that $supp\mu_{Z}\supseteq \partial\sigma(Z)=\left\{\lambda\in
\cc:\,\, |\lambda-1|^2= \frac{ |\alpha\beta||\lambda|}{2}\right\}$.
For this purpose, we need only to prove that $supp\mu_{Z-1}\supseteq
\partial\sigma(Z-1)=\left\{\lambda\in \cc:\,\, |\lambda|^2= \frac{
|\alpha\beta||\lambda+1|}{2}\right\}$.\\
\vskip 0.5cm
 Note that  $\Delta (1+\alpha
E_{12})=1$. For $\lambda\in \cc$, we have \[\log \Delta
((Z-1)-\lambda)=\log \Delta ((1+\alpha E_{12})(1+\beta
F_{12})-(1+\lambda) (1+\alpha E_{12})(1-\alpha E_{12}))\]
\[=\log \Delta (1+\alpha E_{12})+\log \Delta (1+\beta
F_{12}-\lambda+(1+\lambda)\alpha E_{12})=\log\Delta
((1+\lambda)\alpha E_{12}+\beta F_{12}-\lambda).\] By
Example~\ref{E:E12+F12}, $\mu_{(1+\lambda)\alpha E_{12}+\beta
F_{12}}=\mu_{\sqrt{|1+\lambda||\alpha\beta|} b}.$ Hence,
\[\log \Delta ((Z-1)-\lambda)=\log\Delta(\sqrt{|1+\lambda||\alpha\beta|}
b-\lambda)\]
\[
=\log\Delta\left(b-\frac{\lambda}{\sqrt{|1+\lambda||\alpha\beta|}}\right)
-\log\frac{|\lambda|}{\sqrt{|1+\lambda||\alpha\beta|}}+\log|\lambda|.\]
Since $b$ is an R-diagonal operator, this implies that
\begin{equation}\log\Delta\left(b-\frac{|\lambda|}{\sqrt{|1+\lambda||\alpha\beta|}}\right)=
\log\frac{|\lambda|}{\sqrt{|1+\lambda||\alpha\beta|}}-\log|\lambda|+\log
\Delta ((Z-1)-\lambda). \label{Eq:6.1}
\end{equation}
Suppose $\lambda_0\in \partial \sigma(Z-1)$ and $\lambda_0\notin
supp\mu_{Z-1}$. Then there is $\delta>0$ such that
$\mathbb{B}(\lambda_0,\delta)\subset \cc\setminus supp\mu_{Z-1}$.
Now $\log\Delta ((Z-1)-\lambda)$ is a harmonic function on
$\mathbb{B}(\lambda_0,\delta)$. Since $\tau((Z-1)^n)=0$ for all
$n=1,2,\cdots$. By Lemma 4.3 of~\cite{H-L}, for $\lambda\in \cc$
such that $|\lambda|\geq r(Z-1)$, $\log\Delta
((Z-1)-\lambda)=\log|\lambda|$. By the uniqueness of harmonic
functions, we have $\log\Delta ((Z-1)-\lambda)=\log|\lambda|$ for
$\lambda\in \mathbb{B}(\lambda_0,\delta)$. By
equation~(\ref{Eq:6.1}), this implies that
\begin{equation}\log\Delta\left(b-\frac{|\lambda|}{\sqrt{|1+\lambda||\alpha\beta|}}\right)=
\log\frac{|\lambda|}{\sqrt{|1+\lambda||\alpha\beta|}}.\label{Eq:6.2}
\end{equation}
Let $r=\frac{|\lambda|}{\sqrt{|1+\lambda||\alpha\beta|}}$. Then
equation~(\ref{Eq:6.2}) implies that
\begin{equation*}\log\Delta\left(b-r\right)=
\log r
\end{equation*}
for $r\in (s,t)\subseteq [0, \frac{1}{\sqrt{2}}]$. Since $b$ is an
R-diagonal operator, this implies that $\log\Delta(b-z)$ is harmonic
on the annulus with inner radius $s$ and outer radius $t$,
$0<s<t<\frac{1}{\sqrt{2}}$. By Theorem 2.1,
$supp\mu_b=\overline{\mathbb{B}(0,\frac{1}{\sqrt{2}})}$. It is a
contradiction.}
\end{Example}

\section{Hyperinvariant subspaces for operators in $\M$}

\begin{Lemma}\label{L:lemma for center element} For $X\in \M$, if $supp\mu_X=\{\lambda\}$, then
$\tau(X^n)=\lambda^n$ for $n=1,2,\cdots$.
\end{Lemma}
\begin{proof} $\tau(X^n)=\int_{supp\mu_X} z^n d\mu_X(z)=\lambda^n$.
\end{proof}

The converse of Lemma~\ref{L:lemma for center element} is not true.
Since for an R-diagonal operator $X$, we have $\tau(X^n)=0$ for
$n=1,2\cdots$.

\begin{Proposition} Let $X=A+B$, where $A\in M_2(\cc)_{(1)}$ and
$B\in M_2(\cc)_{(2)}$. If $A,B$ are not scalar matrices, then
$supp\mu_X$ contains more than two points.
\end{Proposition}
\begin{proof} Suppose $A,B$ are not scalar matrices.  Since
$A+B=\tau(A)1+\tau(B)1+(A-\tau(A)1)+(B-\tau(B)1)$, to show
$supp\mu_X$ contains more than two points, we need only to show
$\mu_{(A-\tau(A)1)+(B-\tau(B)1)}$ contains more than two points. So
we may assume that $\tau(A)=\tau(B)=0$ and $A,B\neq 0$. Assume that
the spectra of $A$ and $B$ are $\lambda_1,-\lambda_1$ and
$\lambda_2,-\lambda_2$.  If $\tau(A^2)=\tau(B^2)=0$, then $A$ and
$B$ are unitarily equivalent to $\alpha E_{12}$ and $\beta F_{12}$
in $M_2(\cc)_{(1)}$ and $M_2(\cc)_{(2)}$, respectively. Thus
$\mu_X=\mu_{\alpha E_{12}+\beta F_{12}}$. By
Example~\ref{E:E12+F12}, $supp\mu_X$ contains more than two points.
Now suppose
 $\tau(A^2)\neq 0$ or $\tau(B^2)\neq 0$. Without loss of generality,
we assume that $\lambda_1^2=\tau(A^2)\neq 0$. Note that
$\tau(A+B)=0$ and
$\tau((A+B)^2)=\tau(A^2)+\tau(B^2)=\lambda_1^2+\lambda_2^2$. If
$\lambda_1^2+\lambda_2^2\neq 0$, by Lemma~\ref{L:lemma for center
element}, $supp\mu_X$ contains more than two points. Suppose
$\lambda_1^2+\lambda_2^2=0$. Then $\tau(B^2)=\lambda_2^2\neq 0$.
Simple computations show that
\[\tau((A+B)^4)=\tau(A^4)+\tau(B^4)+4\tau(A^2)\tau(B^2)=
\lambda_1^4+\lambda_2^4+4\lambda_1^2\lambda_2^2=(\lambda_1^2+\lambda_2^2)^2+2\lambda_1^2\lambda_2^2=2\lambda_1^2\lambda_2^2\neq
0.\] Note that $\tau(A+B)=0$. By Lemma~\ref{L:lemma for center
element}, $supp\mu_X$ contains more than two points.
\end{proof}

\begin{Proposition}\label{P:brown measure} Let $X=AB$, where $A\in M_2(\cc)_{(1)}$ and
$B\in M_2(\cc)_{(2)}$. If $A, B$ are not scalar matrices, then
$supp\mu_X$ contains more than two points.
\end{Proposition}
\begin{proof} Suppose $A,B$ are not scalar matrices. We
consider the following cases:\\

\noindent Case 1. $\tau(A)=\tau(B)=0$ and $A,B\neq 0$. By Theorem
4.1, $AB(\neq 0)$ is an
R-diagonal operator. So $supp\mu_X$ contains more than two points.\\

\noindent Case 2. $\tau(A)=0, \tau(B)\neq 0$ or $\tau(A)\neq 0,
\tau(B)=0$. Without loss of generality, we assume that $\tau(A)=0$
and $\tau(B)\neq 0$. Then $\tau(AB)=0$ and
$\tau(ABAB)=\tau(A^2)\tau(B)$. If $\tau(A^2)\neq 0$, then
$\tau(ABAB)\neq 0$. By Lemma 10.1, $supp\mu_X$ contains more than
two points. If $\tau(A^2)=0$, then $A$ is unitarily equivalent to
$\alpha E_{12}$ in $M_2(\cc)_{(1)}$. By Lemma 4.5,
$\mu_X=\mu_{\alpha E_{12}(B-\tau(B))}$. Since $\alpha
E_{12}(B-\tau(B)(\neq 0)$ is an R-diagonal operator, $supp\mu_X$
contains
more than two points.\\

\noindent Case 3. $\tau(A)\neq 0$ and $\tau(B)\neq 0$. We may assume
that $\tau(A)=\tau(B)=1$. Let $A=1+A_1$ and $B=1+B_1$. Then
$\tau(A_1)=\tau(B_1)=0$.\\

\noindent Subcase 3.1.  $\tau(A_1^2)\neq 0$ or $\tau(B_1^2)\neq 0$.
We may assume that $\tau(A_1^2)\neq 0$.  Simple computation shows
that $\tau(AB)=1$, $\tau(ABAB)=1+\tau(A_1^2)+\tau(B_1^2)$ and
$\tau((AB)^3)=1+3(\tau(A_1^2)+\tau(B_1^2))+9\tau(A_1^2)\tau(B_1^2).$
If $\tau(A_1^2)+\tau(A_2^2)\neq 0$, then $\tau(ABAB)\neq 1$. By
Lemma~\ref{L:lemma for center element}, $supp\mu_X$ contains more
than two points.  If $\tau(A_1^2)+\tau(A_2^2)=0$, then
$\tau(A_2^2)=-\tau(A_1^2)\neq 0$. So $\tau((AB)^3)\neq 1$. By
Lemma~\ref{L:lemma for center element} again, $supp\mu_X$ contains
more than two points.\\

\noindent Subcase 3.2. $\tau(A_1^2)=\tau(A_2^2)=0$. Then $A_1$ and
$A_2$ are unitarily equivalent to $\alpha E_{12}$ and $\beta F_{12}$
in $M_2(\cc)_{(1)}$ and $M_2(\cc)_{(2)}$, respectively.  So
$\mu_X=\mu_{(1+\alpha E_{12})(1+\beta F_{12})}$.  We may assume
that $A=\left(\begin{array}{cc} 1&\alpha\\
0&1
\end{array}\right)_{(1)}$ and $B=\left(\begin{array}{cc} 1&\beta\\
0&1
\end{array}\right)_{(2)}$. By Example~\ref{E:1+E12 times 1+F12}, $supp\mu_X$ contains more
than two points.
\end{proof}

\begin{Corollary} Let $X=AB$ or $X=A+B$, where $A\in M_2(\cc)_{(1)}$ and
$B\in M_2(\cc)_{(2)}$. If $X\neq \lambda 1$, then $X$ has a
nontrivial hyperinvariant subspace relative to $\M$.
\end{Corollary}
\begin{proof} If $X=A+B$ and $A=\lambda
1$ or $B=\lambda 1$, then $X=\lambda 1+ B$ or $X=\lambda 1+ A$. If
$X$ is not a scalar matrix and $\eta$ is an eigenvalue of $X$, then
$ker(X-\eta 1)$ is a nontrivial hyperinvariant subspace of $X$. If
$X=A+B$ and $A,B\neq \lambda 1$, then $supp\mu_X$ contains more than
two points by Proposition~\ref{P:brown measure}. By~\cite{H-S1}, $X$
has a nontrivial hyperinvariant subspace relative to $\M$. If $X=AB$
and $A=\lambda 1$ or $B=\lambda 1$, then $X=\lambda B$ or $X=\lambda
A$. If $X$ is not a scalar matrix and $\eta$ is an eigenvalue of
$X$, then $ker(X-\eta 1)$ is a nontrivial hyperinvariant subspace of
$X$. If $X=AB$ and $A,B\neq \lambda 1$, then $supp\mu_X$ contains
more than two points. By~\cite{H-S1}, $X$ has a nontrivial
hyperinvariant subspace relative to $\M$.
\end{proof}

{\bf Acknowledgements:}\,The authors want to express their deep
gratitude to professor Eric Nordgren  for valuable discussions. The
authors also thank the referee for some useful  suggestions.

 \vspace{.2in}
\noindent {\em E-mail address: } [Junsheng Fang] jfang\@@cisunix.
unh.edu

\noindent {\em E-mail address: } [Don Hadwin] don\@@math.unh.edu

\noindent {\em E-mail address: } [Xiujuan Ma] mxjsusan@@hebut.edu.cn

\end{document}